\documentclass[11pt]{article} 

\usepackage{enumitem}
\usepackage{float}
\usepackage{multirow}
\usepackage{amsmath,amsfonts,amssymb,graphicx,amsthm,color,natbib}
\usepackage{authblk}
\allowdisplaybreaks
\usepackage[]{geometry}
\usepackage[normalem]{ulem}
\usepackage{hyperref}
\usepackage{booktabs}

\theoremstyle{plain}

\newtheorem{theorem}{Theorem}[section]

\newtheorem{proposition}[theorem]{Proposition}

\newtheorem{remark}{Remark}[section]

\newcommand{\bbR}{\mathbb{R}}
\newcommand{\bbN}{\mathbb{N}}

\begin{document}
	
	\title{Goodness-of-fit test for count distributions\\with finite second moment}
	
\author[1]{{Antonio} {Di Noia}}
\author[1]{{Lucio} {Barabesi}}
\author[1]{{Marzia} {Marcheselli}}
\author[1]{{Caterina} {Pisani}}
\author[2]{{Luca} {Pratelli}}
\affil[1]{Department of Economics and Statistics, University of Siena}
\affil[2]{Italian Naval Academy}
\date{}
\setcounter{Maxaffil}{0}
\renewcommand\Affilfont{\itshape\small}
\maketitle
\let\thefootnote\relax\footnotetext{\emph{Email addresses:} antonio.dinoia55(\sout{at})gmail.com (Antonio Di Noia),  lucio.barabesi(\sout{at})unisi.it (Lucio Barabesi), marzia.marcheselli(\sout{at})unisi.it (Marzia Marcheselli), caterina.pisani(\sout{at})unisi.it (Caterina Pisani), luca\_pratelli(\sout{at})marina.difesa.it (Luca Pratelli).}

	\maketitle
	
\begin{abstract}
A goodness-of-fit test for one-parameter count distributions with finite second moment is proposed. 
The test statistic is derived from the $L^1$ distance of a function of the probability generating function of the model under the null hypothesis and that of the random variable actually generating data, when the latter belongs to a suitable wide class of alternatives. The test statistic has a rather simple form and it is asymptotically normally distributed under the null hypothesis, allowing a straightforward implementation of the test. Moreover, the test is consistent for alternative distributions belonging to the class, but also for all the alternative distributions whose probability of zero is different from that under the null hypothesis. Thus, the use of the test is proposed and investigated also for alternatives not in the class. The finite-sample properties of the test are assessed by means of an extensive simulation study.
\end{abstract}


\section{Introduction}
	Count data naturally arise in many applied disciplines such as actuarial science, medicine, biology and economics, among many others. The Poisson distribution is likely to be the most popular model for such type of data mainly for its simplicity. Nevertheless, observations may exhibit over-dispersion, under-dispersion, zero-inflation or heavy tails, thus precluding the use of the Poisson model as a suitable model. A plethora of count distributions have been introduced that can model these features (e.g., \citealt{johnson2005univariate}). Classical examples are the Negative Binomial for over-dispersion and the zero-inflated Poisson for excesses of zeroes. The Poisson-Tweedie family of distributions, which has been studied by several authors with different parametrization (e.g., \citealt{el2011, barabesi2018tempered, baccini2016random, barabesi2014note}), is able to fit a wide range of mean-variance ratio and tail heaviness.
	Moreover, in order to model over-dispersed data, \cite{tsylova2017using} and \cite{castellares2018bell} introduced the one-parameter Bell family of distributions on the basis of the well-known Bell series expansion (\citealt{bell1934numbers}). These laws have many appealing properties, since they are members of the one-parameter exponential family and are infinitely divisible. A further distribution which has many interesting applications in the setting of queueing theory and branching processes (\citealt{johnson2005univariate}) is the Borel law (\citealt{borel1942emploi}).
	
	A challenging aspect of data analysis consists in testing the goodness-of-fit to a parametric family of count distributions. Many testing procedures dealing with count data are based on the properties of the probability generating function (p.g.f.) and on the corresponding empirical p.g.f. Indeed, the p.g.f. fully characterizes the distribution, it is sometimes simpler than the corresponding probability mass function (p.m.f.) and possesses convenient features, since it is a real-valued continuous function always defined in the range $[0,1]$. The use of the p.g.f. in testing the fit of discrete distributions has a long-standing tradition (e.g., \citealt{kocherlakota1986goodness, rueda1991goodness}). In particular, \cite{rueda1991goodness} introduced a test for the Poisson distribution with known parameter, extended by \cite{rueda1999tests} to the case of unknown parameter and to the negative Binomial distribution. As to testing Poissonity, \cite{nakamura1993use} proposed a test based on the empirical p.g.f., while \cite{meintanis2008class} and \cite{puig2020some} introduced tests based on different characterizations of the p.g.f. against alternatives belonging to a large
	family. In a more general framework, \cite{jimenez2017minimum} presented a test statistic based on a distance between the empirical p.g.f. and the p.g.f. of the model under the null hypothesis, together with a weighted bootstrap estimator of its distribution. Moreover, \cite{jimenez2019testing} introduced a computationally convenient test for the Poisson–Tweedie distribution, while \cite{jimenez2021test} suggested a test for the Geometric distribution.
	
	In this paper, a novel goodness-of-fit test for families of one-parameter count distributions with finite second moment is proposed.
	The test stands in the long tradition of testing procedures based on distances, such as the Pearson chi-squared test. In particular, the proposed test statistic is justified by the $L^1$ distance of a suitable function of the p.g.f. of the model under the null hypothesis and the p.g.f. of the random variable actually generating data, when the derivative of the ratio of the p.g.f.s has constant sign. Therefore, given the distribution specified under the null hypothesis, the natural class of alternative distributions contains those ensuring the derivative constant sign and, for the corresponding hypothesis system, the test is proven to be consistent. 
    The test statistic has a manageable expression and depends on the empirical p.g.f. solely through its value in zero, thus avoiding the complexities of handling the whole empirical functional. In addition, the test statistic is proven to have an asymptotic normal distribution, which allows for a straightforward implementation of the test, without demanding intensive resampling methods. Moreover, the test can be also adopted for the very general hypothesis system with alternatives not necessarily belonging to the class, even though in this case the consistency is ensured only if the probability of getting zero is different under the null and alternative hypothesis.  
     	
	Section \ref{sec:prelim} contains some preliminaries about the hypothesis system and the distance criterion. In Section \ref{sec:test-stat}, the new goodness-of-fit test is proposed and its asymptotic properties are proven. Section \ref{sec:families} deals with the test statistics for some well-known families of count distributions. In Section \ref{sec:contig} the asymptotic behaviour of the proposed test is investigated under contiguous alternatives. A Monte Carlo simulation to assess the finite-sample performance of the test is described in Section \ref{sec:sim}. Some concluding remarks are given in Section \ref{sec:disc}. All tables and figures are reported in the Appendix.

	\section{Preliminaries}
	\label{sec:prelim}
	Let $\Theta$ be a subset of $\bbR$, $\{{\cal M}_\theta\}_{\theta\in\Theta}$ a family of distributions concentrated into ${\bbN}_0$ with finite second moment and $p_{\theta}$ the corresponding p.m.f. Without loss of generality, let assume ${\cal M}_\theta(\{0\})\neq 0$, i.e. the singleton $\{0\}$ is not negligible with respect to ${\cal M}_\theta$. 
    Examples of widely applied families of type ${\cal M}_\theta$ are the Poisson family $\{{\cal P}(\lambda)\}_{\lambda\in\bbR^+}$, the Geometric family $\{{\cal G}(p)\}_{p\in]0,1[}$, the Bell family $\{{\cal B}e(\theta)\}_{\theta>0}$ and the shifted Borel family $\{{\cal B}o(\lambda)\}_{\lambda>0}.$ 
    
    When interest is in assessing if the random variable (r.v.) $X$ is distributed according to the model $\{{\cal M}_\theta\}_{\theta\in\Theta}$, i.e. $H_0: X \sim {\cal M}_\theta$ for some $\theta\in\Theta$, a convenient parametrization may be achieved by using the first moment $\mu$, in order to derive the asymptotic properties of suitable test statistics more easily by means of limit theorems.
    In particular, if $h: \theta \mapsto \sum_{n=0}^\infty n p_\theta(n)$
	is a strictly monotone function, the model is parametrized as ${\cal M}_{h^{-1}(\mu)}$.
	The Geometric, the shifted Borel and the Bell family can be parametrized through $\mu$ since $p=(\mu+1)^{-1}$, $\lambda={{\mu}\over{1+\mu}}$ and $\theta=h^{-1}(\mu)$, where $h^{-1}$ is the inverse function of  $\theta\mapsto \theta e^\theta$. 
	
	Let $f$ represent the p.g.f. of the r.v. $X$ and $g_\mu$ be the p.g.f. under the model ${\mathcal M_{h^{-1}(\mu)}}$. There exist various proposals in literature (e.g., \citealt{sim2010parameter} and references therein) for quantifying discrepancy between $f$ and $g_\mu$  and a further sensible measure could be based on the ratio $f/g_\mu$ or on the corresponding derivative $(f/g_\mu)^\prime$. Indeed, since  $$\Big(\frac{f}{g_\mu}\Big)^\prime=\frac{f}{g_\mu}\Big(\frac{f^\prime}{f}-\frac{g_\mu^\prime}{g_\mu}\Big),$$
	the derivative could be considered as a \lq\lq weighted\rq\rq version of the original ratio, where the weight is given by the difference of the normalized variation of the single p.g.f.s. In literature the normalized variation is used to give an interesting characterization of infinitely divisible p.g.f.s (see Theorem 4.2 in \citealt{steutel2003infinite}). Moreover, the difference of the normalized variation (with the appropriate sign) could be more effective to detect small discrepancies between $f$ and $g_\mu$. 
	Obviously, if $X$ is distributed according to ${\mathcal M_{h^{-1}(\mu)}}$, then $(f/g_\mu)(s)=1$ and $(f/g_\mu)^\prime(s)=0$ for any $s \in [0,1]$. 
	Hence, denoting by
	$$D_\mu(s)=\Big({{f}\over{g_{\mu}}}\Big)^\prime(s),$$
	the $L^1([0,1])$ distance of $D_\mu(s)$ from the null function can be considered, since such a distance is zero under $H_0$ while positive values evidence departures from $H_0$. Assuming that $f$ and $g_{\mu}$ are such that $D_\mu(s)$ is non-negative or non-positive for any $s\in[0,1]$, the $L^1([0,1])$ distance is defined as
	
	\begin{equation}
		\label{ds}
		\int_0^1 |D_\mu(s)|ds =\Big\vert\int_0^1 D_\mu(s) ds\Big\vert=\Big\vert{{f}\over{g_{\mu}}}(1)-{{f}\over{g_{\mu}}}(0)\Big\vert   =\Big\vert{{g_{\mu}(0)-P(X=0)}\over{g_{\mu}(0)}}\Big\vert
	\end{equation}
    and thus a reasonable test statistic for assessing $H_0$ could be based on suitable estimators of $\mu$ and $P(X=0)$.  
	
	From \eqref{ds}, it is also natural to consider the class $\Delta_{{\mathcal M}_{h^{-1}(\mu)}}$ of count distributions (which are not in  $H_0$), depending on the model specified under $H_0$, such that $D_\mu(s)$ has constant sign for any $s\in[0,1]$, and the corresponding hypothesis system
   $H_0: X \sim {\mathcal M}_{h^{-1}(\mu)} ,\ H_1: X\sim \Delta_{{\mathcal M}_{h^{-1}(\mu)}}$. Considering a fairly wide
   class of alternatives has already been exploited in literature (see e.g. \citealt{meintanis2008class}, \citealt{puig2020some}). Also in this setting the class $\Delta_{{\mathcal M}_{h^{-1}(\mu)}}$ is rather wide: 
   for example, for the geometric family the class contains many widely applied distributions such as the Poisson, the Binomial, the Negative Binomial and the Neyman type A distribution.  

	It must be pointed out that previous hypothesis system arises in many contexts where the assessment of the null hypothesis may be difficult. In particular, 
	$H_0$ may be hard to assess when $f$ is equal to $g_{\mu_0}w$ and $w$ is in turn a p.g.f. in such a way that
	 $f$  constitutes a \lq\lq perturbation\rq\rq of $g_{\mu_0}$, i.e. $f$ nearly resembles $g_{\mu_0}$ when $w(0)$ is close to one. 
	Obviously, $D_{\mu_0}(s)$ is non-negative for any $s \in [0,1]$ and $X$ is distributed as $X_0\sim {\cal M}_{h^{-1}(\mu_0)}$ under $H_0$, while $X=X_0+Y$ under $H_1$, where $Y$ is a r.v.  with p.g.f. $w$, independent of $X_0$. An interesting case is obtained when $Y=I_{\{Z\geq 1\}}\sum_{n=1}^{Z}X_n$, i.e. $X_0+Y$ is the random sum $\sum_{n=0}^{Z}X_n$ where $(X_n)_n$ is a sequence of independent r.v.s with $X_n\sim X_0$ and $Z$ is a non-negative integer-valued r.v. independent of $X_n$  for any natural $n$.

    The null hypothesis can be  difficult to assess also when $f$ is the p.g.f. of the r.v. $\alpha$-fraction of $X_0$ given by $\sum_{n=1}^{X_0} Y_n$, where $(Y_n)_n$ are i.i.d. Bernoulli r.v.s with parameter $\alpha\in[0,1]$, independent of $X_0$. Following \cite{steutel2003infinite}, the $\alpha$-fraction of $X_0$ is defined by means of the so-called binomial thinning operator. It is worth noting that $g_{\mu_0}\leq f$ and $\sum_{n=1}^{X_0} Y_n$ converges almost surely to $X_0$ for $\alpha$ approaching 1. Furthermore, $D_{\mu_0}(s)$ is non-positive for any $s \in [0,1]$ for $X_0$ belonging to many families of type  ${\mathcal M_{h^{-1}(\mu_0)}}$, such as those of the Binomial, the Negative Binomial, the Logarithmic, the Sibuya, the discrete stable, the discrete Linnik.

	\section{The test statistics}
	\label{sec:test-stat}
	 Thanks to \eqref{ds}, we can introduce a family of test statistics based on estimators of $\mu$ and $P(X=0)$. To this aim, given a random sample $X_1,\ldots,X_n$ from $X$ and denoting by $\widehat{\mu}_n$ and $\widehat{P}_n(0)$ the sample mean and the sample proportion of observations equal to zero, that is 
	
	\begin{equation*}
		\widehat\mu_n=\frac{X_1+\ldots+X_n}{n},\  \widehat{P}_n(0)=\frac{I_{\{X_1=0\}}+\ldots+I_{\{X_n=0\}}}{n},
	\end{equation*}
	a test statistic can be based on  
	
	\begin{equation*}
		\label{T0_hat}
		\widehat{T}_0=\sqrt n(g_{\widehat\mu_n}(0)-\widehat{P}_n(0)).
	\end{equation*}

	In the following proposition we prove that the asymptotic distribution of $\widehat{T}_0$ under the null hypothesis does not rely on peculiar characteristics of the class $\Delta_{{\mathcal M}_{h^{-1}(\mu)}}$. Then, we propose its use for assessing the much more general hypothesis system $$H_0: X \sim {\mathcal M}_{h^{-1}(\mu)}\,\, \text{for some}\,\, {\mu\in {\bbR^+}}, \ H_1: X \nsim {\mathcal M}_{h^{-1}(\mu)} \,\, \text{for all}\,\, {\mu\in {\bbR^+}},$$ even if its interpretation as deriving from a $L^1$ distance is lost. Referring to this more general hypothesis system, $\widehat{T}_0$ can be rewritten as
$$\widehat{T}_0=\sqrt n(\psi_0 (\widehat\mu_n)-\widehat{P}_n(0)),$$
where $\psi_0:\mu\mapsto {\cal M}_{h^{-1}(\mu)}(\{0\})$ and its asymptotic distribution is derived under some mild conditions on $\psi_0$.

	\begin{proposition}
		\label{proposition1}
		Let $\psi_0$ be a $C^1$ function with bounded first-order derivative. Then, under the null hypothesis,  
		$\widehat{T}_0$ converges in distribution to ${\mathcal N}(0,\sigma^2)$ as $n\to \infty$, where $\mu=E[X]$ and
		\begin{equation}
			\label{sigma2_1}
			\begin{aligned}
				\sigma^2&={\rm Var}[\psi_0^\prime(\mu)X-I_{\{X=0\}}]\cr
				&=(\psi_0^\prime)^2(\mu){\rm Var}[X]+2\psi_0^\prime(\mu)\mu\psi_0(\mu)+\psi_0(\mu)(1-\psi_0(\mu)).
			\end{aligned}
		\end{equation}
		Moreover, if $X \sim \Delta_{{\mathcal M}_{h^{-1}(\mu)}}$ or, more in general, if $r_0=\psi_0(\mu)-P(X=0)\not=0$, then $|\widehat{T}_0|$ converges in probability to $\infty$.
	\end{proposition}
	
	\begin{proof}
		Owing to the Delta Method 
		\begin{equation*}
			\begin{aligned}
				\widehat{T}_0&= \sqrt n(\psi_0(\widehat\mu_n)-\psi_0(\mu)+\psi_0(\mu) -\widehat{P}_n(0))\cr 
				&=\sqrt n(\psi_0^\prime(\mu)(\widehat\mu_n-\mu)+\psi_0(\mu)-\widehat{P}_n(0))+o_P(1)\cr 
				&={\frac{\tau(X_1)+\ldots+\tau(X_n)}{\sqrt n}}-\sqrt n r_0+o_P(1),
			\end{aligned}
		\end{equation*}
		where $\tau$ is the function defined by $$x\mapsto \psi_0^\prime(\mu)(x-\mu)-(I_{\{X=0\}}-P(X=0)).$$
		Under $H_0$, $r_0=0$, ${\rm E}[\tau(X)]=0$ and ${\rm E}[o_P^2(1)]$ is $o(1)$ since $$|\psi_0(\widehat\mu_n)-\psi_0(\mu)|\leq |\widehat\mu_n-\mu|.$$ Then, under $H_0$, 
		$\widehat{T}_0$ converges in distribution to ${\mathcal N}(0,{ {\rm Var}\, [\tau(X)]})$ by applying the Central Limit Theorem to ${\frac{\tau(X_1)+\ldots+\tau(X_n)}{\sqrt n}}$. Moreover, since
		\begin{equation*}
			\begin{aligned}
				{\rm Var}[\tau(X)]&={\rm Var}[\psi_0^\prime(\mu)X-I_{\{X=0\}}]\cr
				&=(\psi_0^\prime)^2(\mu){\rm Var}[X]+{\rm Var}[I_{\{X=0\}}]-2\psi_0^\prime(\mu){\rm Cov}[X,I_{\{X=0\}}]	
			\end{aligned}
		\end{equation*}
		and
		$${\rm Cov}[X,I_{\{X=0\}}]=-\mu \psi_0(\mu)$$ the first part of the proposition is proven.
		
		Now, let $X$ be a r.v. such that $r_0\neq 0$. Since $${\frac{\tau(X_1)+\ldots+\tau(X_n)}{\sqrt n}}+o_P(1)$$ is bounded in probability and $\sqrt n |r_0|$ converges to $\infty$, then
		 $|\widehat{T}_0|$ converges in probability to $\infty$. The second part of the proposition is so proven.
	\end{proof}
	
	It is worth noting that $\sigma^2>0$ when $\psi_0^\prime(\mu)=0$. Moreover, if $\psi_0^\prime(\mu)\not=0$, $\sigma^2=0$ iff there exists a real number $c$ such that $\psi_0^\prime(\mu)X-I_{\{X=0\}}=c$ almost surely. Since $P(X=0)$ is not negligible, $\sigma^2=0$ iff $X=-I_{\{X\not=0\}}/\psi_0^\prime(\mu)$. Thus, $\sigma^2>0$ if $X$ takes more than two values, as it happens for all the distributions considered in the following sections, and therefore it is not restrictive to consider $\sigma^2>0$.

	In order to obtain a test statistic, $\sigma^2$ can be estimated by means of the plug-in estimator
	
	\begin{equation}
		\label{sigma2_1_hat}
		\widehat\sigma^2_n=(\psi_0^\prime)^2(\widehat\mu_n)
		v(\widehat\mu_n)+2\widehat\mu_n\psi_0(\widehat\mu_n)
		\psi_0^\prime(\widehat\mu_n)+\psi_0(\widehat\mu_n)(1-\psi_0(\widehat\mu_n))
	\end{equation}
	where
	\begin{equation*}
		v(\widehat\mu_n)=\sum_{j=0}^\infty (j-\widehat\mu_n)^2 p_{h^{-1}(\widehat\mu_n)}(j).
	\end{equation*} 
	Since $\widehat\sigma^2_n$ converges almost surely to $\sigma^2$ it follows $\lim_n \widehat\sigma^2_n>0$ almost surely. Thus the test statistic is defined as
	\begin{equation*}
		Z_{n}=\frac{\widehat T_0}{\widehat\sigma_{n}}
	\end{equation*}
    if $\widehat\sigma_{n}>0$ and $Z_{n}=0$ else, which is asymptotically equivalent to $\frac{\widehat T_0}{\sigma}$ and, thanks to Proposition \ref{proposition1}, it has an asymptotic ${\mathcal N}(0,1)$ distribution. It is at once apparent that the rejection region of an $\alpha$-level large-sample test is given by
    $\{|Z_{n}|>z_{ 1-\alpha/2}\},$ 
    where $z_{ 1-\alpha/2}$ is the $({1-\alpha/2})$-quantile of the standard normal distribution.   
	Moreover, the test is consistent for the alternatives in $\Delta_{{\mathcal M}_{h^{-1}(\mu)}}$ but also for all the other alternatives not in the class for which the probability of $0$ is different from that under the null hypothesis. Requiring the probability of $0$ being different may be a rather restrictive assumption even though the class of families of distributions not in $\Delta_{{\mathcal M}_{h^{-1}(\mu)}}$ may be narrow.

	Obviously, alternative suitable estimators of $\sigma^2$ are possible for particular families of distributions or by using maximum likelihood estimators of $\mu$.

	\section{The test statistic for some families of distributions}
	\label{sec:families}
	\subsection{Shifted Borel family}
	The Borel family arises in the context of queueing theory and branching processes. More precisely, the Borel distribution (\citealt{borel1942emploi}) describes the distribution of the total number of customers served before a queue vanishes, given a single queue with Poisson random arrival of customers and a constant time in serving each customer, when there is initially one customer in the queue. The Borel distribution has parameter $\lambda$ if the constant time is $1$ and the constant rate of arrivals is $\lambda$. Equivalently, the Borel distribution is the distribution of the total progeny of a Galton-–Watson branching
	process where each individual has ${\mathcal P}(\lambda)$ children (e.g., \citealt{borel1942emploi, janson2008susceptibility, johnson2005univariate}). In particular, these distributions are concentrated on ${\bbN}_0$ when $\lambda\leq 1$ and have finite moment of any order when $\lambda<1$.

	In the following, the family of shifted Borel distributions $\{{\mathcal B}o(\lambda)\}_{\lambda\in [0,1[}$ with values in ${\bbN_0}$ and p.m.f.
	$$P(X=n)=e^{-\lambda (n+1)}{{(\lambda(n+1))^n}\over{(n+1)!}} \ \ \ \ \ n\in {\bbN_0}$$ 
	is considered.
	Since $\mu={{\lambda}\over{1-\lambda}}=h(\lambda)$, then $h^{-1}(\mu)={{\mu}\over{1+\mu}}$. Moreover, $P(X=0)=\psi_0(\mu)=e^{-{{\mu}\over{1+\mu}}}$ and
	$${\widehat T}_0=\sqrt n\big(e^{-{{\widehat\mu_n}\over{1+\widehat\mu_n}}}-
	\widehat{P}_n(0)\big).$$ 
	Owing to Proposition \ref{proposition1}, $\widehat{T}_0$ converges in distribution to ${\mathcal N}(0,\sigma^2)$ as $n\to \infty$, where $$\sigma^2={\rm Var}\Big[\frac{e^{-\frac{\mu}{\mu+1}}}{(\mu+1)^2}X+I_{\{X=0\}}\Big]$$
	and, since ${\rm Var}[X]={{\lambda}\over{(1-\lambda)^3}}=\mu(1+\mu)^2$, from \eqref{sigma2_1} $$\sigma^2=e^{-{{2\mu}\over{1+\mu}}}\Big(e^{{{\mu}\over {1+\mu}}}-1-{{{\mu}\over{(1+\mu)^2}}}\Big) $$ and from \eqref{sigma2_1_hat} $$\widehat\sigma^2_n=e^{-{{2\widehat\mu_n}\over{1+\widehat\mu_n}}}\Big(e^{{{\widehat\mu_n}\over {1+\widehat\mu_n}}}-1-{{{\widehat\mu_n}\over{(1+\widehat\mu_n)^2}}}\Big).$$
	
	\subsection{Geometric family}
	Let us consider the well-known family of Geometric distributions $\{{\mathcal G}(p)\}_{p\in]0,1]}$. The r.v. $X$ has p.m.f. given by
	$$P(X=n)=p(1-p)^n\ \ \ \ \ n\in {\bbN_0}.$$ 
	Since $\mu={{1-p}\over{p}}=h(p)$, it holds $h^{-1}(\mu)=\frac{1}{\mu+1}$ and $P(X=0)=\psi_0(\mu)=\frac{1}{\mu+1}.$ Then 
	$${\widehat T}_0=\sqrt n\Big(\frac{1}{\widehat\mu_n+1}-
	\widehat{P}_n(0)\Big)$$ 
	and $\widehat{T}_0$ converges in distribution to ${\mathcal N}(0,\sigma^2)$ as $n\to \infty$, where $$\sigma^2={\rm Var}\Big[\frac{X}{(\mu+1)^2}+I_{\{X=0\}}\Big].$$
	From Proposition \ref{proposition1}, since ${\rm Var}[X]={{1-p}\over{p^2}}={{\mu}{(1+\mu)}}$, from \eqref{sigma2_1} and from \eqref{sigma2_1_hat} it follows
	$$\sigma^2=\frac{\mu^2}{(\mu+1)^3}$$
	and
	$$\widehat\sigma^2_n=\frac{\widehat\mu_n^2}{(\widehat\mu_n+1)^3}.$$
	
	\subsection{Bell family}
	The Bell family of distributions $\{{\mathcal B}e(\theta)\}_{\theta>0}$  has been recently introduced by  \cite{tsylova2017using} and \cite{castellares2018bell}. Bell distributions have many interesting properties. Indeed, the Bell family belongs to the exponential family and it is infinitely divisible. Moreover, the Poisson distribution is not nested in the Bell family but it can be approximated for small values of the parameter by the Bell distribution. This family is rather flexible for fitting a wide spectrum of count data which may present over-dispersion and it may be an alternative model to the very familiar Poisson and Negative Binomial models in several areas. For example, owing to its flexibility, the Bell distribution is used to model the number of insurance claims over a fixed period of time. For further applications see \cite{batsidis2020goodness}.
	
	The Bell p.m.f. with parameter $\theta$ has a very simple form and it is given by
	$$P(X=n)={{{\theta^nB_ne^{1-e^{\theta}}}}\over{n!}}\ \ \ \ \ n\in {\bbN_0},$$
	where the Bell number $B_n$ (see \citealt{bell1934numbers}) is the number of partitions of a set of size $n$ and is equal to the $n$-th moment of a Poisson r.v. with $\mu=1$.
	Since $\mu=\theta e^\theta=h(\theta)$, $h^{-1}$ does not have closed form but simple numerical techniques can be adopted to obtain the value of the inverse function at $\mu$ or at its estimate. Thus  
	$P(X=0)=\psi_0(\mu)=e^{1-e^{h^{-1}(\mu)}}$
	and
	$${\widehat T}_0=\sqrt n\big(e^{1-e^{h^{-1}(\widehat\mu_n)}}-
	\widehat{P}_n(0)\big).$$ 
	From Proposition \ref{proposition1}, $\widehat{T}_0$ converges in distribution to ${\mathcal N}(0,\sigma^2)$ as $n\to \infty$, where $$\sigma^2={\rm Var}\Big[\frac{e^{1-e^{h^{-1}(\mu)}}}{(1+h^{-1}(\mu))}X+I_{\{X=0\}}\Big].$$
	Moreover, since ${\rm Var}[X]=\theta e^\theta (1+\theta)=\mu(1+h^{-1}(\mu))$, from \eqref{sigma2_1} and from \eqref{sigma2_1_hat} it follows
	$$\sigma^2=e^{1-e^{h^{-1}(\mu)}}\Big(1-{{h^{-1}(\mu)e^{1+h^{-1}(\mu)-e^{h^{-1}(\mu)}}}\over {1+h^{-1}(\mu)}}-e^{1-e^{h^{-1}(\mu)}}\Big)$$
	and
	$$\widehat\sigma^2_n=e^{1-e^{h^{-1}(\widehat\mu_n)}}\Big(1-{{h^{-1}(\widehat\mu_n)e^{1+h^{-1}(\widehat\mu_n)-e^{h^{-1}(\widehat\mu_n)}}}\over {1+h^{-1}(\widehat\mu_n)}}-e^{1-e^{h^{-1}(\widehat\mu_n)}}\Big).$$
	
	\section{Asymptotic behaviour under contiguous alternatives}
	\label{sec:contig}
	
	The asymptotic behaviour of the test statistic $Z_n$ is investigated under suitable contiguous alternatives, obtained by mixtures of distributions. The concept of contiguity is frequently applied in many asymptotic settings (e.g., \citealt{van2000asymptotic, dhar2016study, betsch2019new, kalemkerian2020independence, meselidis2020statistical} among others). In par\-ticular,  let $\{A_{l,n}\}_{l=1,\ldots,n}$ be a triangular array of independent events and  $(Y_n)_n$ be a sequence of i.i.d. non-negative integer-valued r.v.s with $ {\rm E}[Y^{2}_1]<\infty$. Moreover, suppose $(I_{A_{l,n}})_{l=1,\ldots,n}$ and  $Y_1,\ldots, Y_n$ to be mutually independent and also independent of the i.i.d. random variables $X_1,\ldots,X_n$, where now $X_1\sim {\mathcal M_{h^{-1}(\mu)}}$ with $\mu={\rm E}[X_1]$ and $X_1$ takes more than two values. For $l=1,\ldots,n$, denote by
	\begin{equation}
		\label{Xprime} 
	X^\prime_{l,n}=I_{A_{l,n}}X_l+I_{A_{l,n}^c}Y_l  
	\end{equation}
	\noindent with $P(A_{l,n}^c)={{\lambda}\over{\sqrt n}}>0$. Given the random sample $X^\prime_{1,n},\ldots, X^\prime_{n,n}$, let $$\widetilde{\mu}_n={{X^\prime_{1,n}+\ldots+X^\prime_{n,n}}\over{n}},\qquad  \widetilde{P}_n(0)={{I_{\{X^\prime_{1,n}=0\}}+\ldots+I_{\{X^\prime_{n,n}=0\}}}\over{n}}$$ be the corresponding sample mean and sample proportion.  
	It is at once apparent that also $X^\prime_{l,n}$ is a non-negative integer-valued r.v. which converges in $L^2$ to $X_l$ for $n$ approaching infinity. 

	The following result is useful to highlight the discriminatory capability of the test statistic under non-trivial contiguous alternatives.
	
	\begin{proposition}
		\label{proposition2} 
		Let $\psi_0:\mu\mapsto {\cal M}_{h^{-1}(\mu)}(\{0\})$ be a $C^1$ function with bounded first-order derivative and $$\widetilde\sigma^2_n=(\psi_0^\prime)^2(\widetilde\mu_n)
		v(\widetilde\mu_n)-2\widetilde\mu_n\psi_0(\widetilde\mu_n)
		\psi_0^\prime(\widetilde\mu_n)+\psi_0(\widetilde\mu_n)(1-\psi_0(\widetilde\mu_n)).$$
		If 
		${\rm E}[Y_1]={\rm E}[X_1]$ then $\widetilde\sigma^2_n$ converges in probability to ${\rm Var}[\psi_0^\prime({\mu})X_1-I_{\{X_1=0\}}]>0$ and
		\begin{equation}
			\label{tildeT_0}
		{{\widetilde{T}_0-\lambda\big(P(X_1=0)-P(Y_1=0)\big)}\over{\widetilde\sigma_n}} I_{\{\widetilde\sigma_n>0\}}	\stackrel{d}{\longrightarrow} {\mathcal N}(0,1)
		\end{equation}
	 where
		$$\widetilde{T}_0=\sqrt n(\psi_0(\widetilde\mu_n)-\widetilde{P}_n(0)).$$

	\end{proposition}
	
	\begin{proof}
		Consider the difference between $\widetilde{T}_0-\lambda\big(P(X_1=0)-P(Y_1=0)\big)$ and $\widehat{T}_0$, that is
		\begin{equation*}
			R_n=\widetilde{T}_0-\widehat{T}_0-\lambda\big(P(X_1=0)-P(Y_1=0)\big).
		\end{equation*}
		Note that 
			\begin{equation*}
					\begin{aligned}
			R_n&=\sqrt n\big(\psi_0(\widetilde{\mu}_n)-\widetilde{P}_n(0)\big)-\sqrt n\big(\psi_0(\widehat{\mu}_n)-\widehat{P}_n(0)\big)-\lambda\big(P(X_1=0)-P(Y_1=0)\big)\cr & = \psi_0^\prime(\widehat\mu_n)(\widetilde\mu_n-\widehat\mu_n)-(\widetilde{P}_n(0)-\widehat{P}_n(0))-\lambda\big(P(X_1=0)-P(Y_1=0)\big)+o_P(1).
					\end{aligned}	
					\end{equation*}
	   Since
		$$\widetilde\mu_n-\widehat\mu_n={{\sum_{l=1}^nI_{A_{l,n}^c} (Y_l-X_l)}\over{n}}$$
		and 
		$$\widetilde{P}_n(0)-\widehat{P}_n(0)={{\sum_{l=1}^n I_{A_{l,n}^c} (I_{\{Y_{l}=0\}}-I_{\{X_{l}=0\}})}\over{n}},$$
		then
		\begin{equation*}
		\begin{aligned}
			R_n=\psi_0^\prime(\widehat\mu_n){{\sum_{l=1}^n\kern-0.5mm I_{A_{l,n}^c} \kern-0.5mm(Y_l\kern-0.5mm-\kern-0.5mm X_l)}\over{\sqrt n}}-{{\sum_{l=1}^n \big(\kern-0.5mm I_{A_{l,n}^c}\kern-0.5mm(I_{\{Y_{l}=0\}}\kern-0.5mm-\kern-0.5mm I_{\{X_{l}=0\}})-a_n\big)}\over{\sqrt n}}+o_P(1),\cr
		\end{aligned}
	   \end{equation*}
     where $a_n=\lambda\big(P(Y_1=0)-P(X_1=0)\big)/\sqrt n = {\rm E}[I_{A_{l,n}^c}(I_{\{Y_{l}=0\}}\kern-0.5mm-\kern-0.5mm I_{\{X_{l}=0\}})].$
     
     \noindent Moreover, as $${\rm E}\Big[\Big({{\sum_{l=1}^n\kern-0.5mm I_{A_{l,n}^c} \kern-0.5mm(Y_l\kern-0.5mm-\kern-0.5mm X_l)}\over{\sqrt n}}\Big)^2\Big]={{\lambda {\rm E}[(Y_1\kern-0.5mm-\kern-0.5mm X_1)^2]}\over{\sqrt n}}$$ and 	
		\begin{equation*}
		\begin{aligned}
			{\rm E}\Big[\Big({{\sum_{l=1}^n I_{A_{l,n}^c}(I_{\{Y_{l}=0\}}\kern-0.5mm-\kern-0.5mm I_{\{X_{l}=0\}})-a_n}\over{\sqrt n}}\Big)^2\Big] &={\rm Var}\big[I_{A_{1,n}^c}\kern-0.5mm(I_{\{Y_{1}=0\}}\kern-0.5mm-\kern-0.5mm I_{\{X_{1}=0\}})\big]\cr &\leq {{\lambda}\over{\sqrt n}},\cr
		\end{aligned}
	\end{equation*}
    $R_n$ converges in probability to $0$. Therefore, also $R_n/\widehat\sigma_n$ converges in probability to $0$ as $\widehat\sigma_n$ converges in probability to ${\rm Var}[\psi_0^\prime({\mu})X_1-I_{\{X_1=0\}}]$ which is positive and 
    $${{{\widetilde T}_0-\lambda\big(P(X_1=0)-P(Y_1=0)\big)}\over{\widehat\sigma_n}}$$ has the same asymptotic distribution of  ${\widehat T}_0/\widehat\sigma_n$. As $\widetilde\sigma_n$ and $\widehat\sigma_n$ are asymptotically equivalent, since they both converge in probability to  ${\rm Var}[\psi_0^\prime({\mu})X_1-I_{\{X_1=0\}}]$, the convergence of \eqref{tildeT_0} is proven.
    \end{proof}

\begin{remark}Given a family of i.i.d. Bernoulli r.v.s $(Y_{j,l})_{j,l}$ with parameter $1-{{\lambda}\over{\sqrt n}}$ independent of $(X_n)_n$, let $X^\prime_{l,n}=\sum_{j=1}^{X_l} Y_{j,l}$ be 
the $\alpha$-fraction of $X_l$, with $l=1,\ldots,n$. Since ${\rm E}[X^\prime_{l,n}]-\mu$ is equal to $-\mu\lambda/{\sqrt n}$, following the same reasoning of Proposition \ref{proposition2}, 
$${{\widetilde{T}_0+\lambda\big(P(X_1=1)+\widetilde\mu_n\big)}\over{\widetilde\sigma_n}}I_{\{\widetilde\sigma_n>0\}}\stackrel{d}{\longrightarrow} {\mathcal N}(0,1).$$
\end{remark}

	\begin{remark} In both propositions the asymptotic behaviour still holds also removing the condition of bounded first-order derivates, even though the order of convergence of $\widehat\sigma_n$ could be considerably reduced. 
    \end{remark} 
	
	\section{Simulation study}
	\label{sec:sim}
	The performance of the proposed test has been assessed and compared to that of the chi-squared goodness-of-fit test, by means of an extensive Monte Carlo simulation, when the Geometric, Bell and Borel distributions are specified under the null hypothesis. The chi-squared test is suitable for the general hypothesis system and, similarly to our proposal, it is based on a test statistic having known asymptotic distribution and not requiring intensive resampling methods. It is worth noting that the comparison is also meaningful since both tests rely on distance-based statistics. Moreover, we compared the performance of the test with that of a recently introduced test based on bootstrap procedures and specifically tailored for the Geometric  distribution (\citealt{jimenez2021test}). 
		
	As to the chi-squared test, it is well-known that the asymptotic approximation is usually satisfactory if each expected frequency is large enough, with some authors suggesting the minimum value of $5$ and some others of $3$. Therefore, there is no way of avoiding the arbitrariness of grouping classes to compute the chi-squared statistic (e.g., \citealt{gibbons2020nonparametric}). Moreover, the standard implementation considers the sample maximum and the sample minimum as extreme values for the test statistic computation and thanks to simulation studies it is well-known that this choice is not enough in order to capture a sufficient probability mass under $H_0$.  
    Hence, to ensure a more reliable approximation, an ad-hoc general version of the chi-squared test is implemented. In particular, denoting by $\{C_1,\dots,C_k\}$ the set of classes to be considered, $C_1$ is the set of natural numbers smaller than the largest integer not greater than $\mu - 3\sqrt{\mu}$, $C_k$ is the set of natural numbers greater than the smallest integer not less than $\mu + 3\sqrt{\mu}$ and $\{C_2,\dots, C_{k-1}\}$ are the singletons not included in $C_1$ and $C_k$. 
    Obviously, in order to implement the chi-squared test statistic, the classes ${\widehat C}_1,\ldots,{\widehat C}_k$ and the corresponding expected frequencies $e_1,\ldots,e_k$  are obtained by plugging the parameter estimates, also adopted in $Z_n$, in the null distribution. In particular, the chi-squared test statistic $Q_n$ is given by
    \begin{equation*}
    Q_n=\sum_{j=1}^k\frac{(n_j-e_j)^2}{e_j}
    \end{equation*} 
    where $n_j$ denotes the number of observations in the class ${\widehat C}_j$ and $e_j=n\sum_{l\in{\widehat C}_j}p_{h^{-1}(\widehat\mu_n)}(l).$ The simulation is implemented by using R (\citealt{Rcore}) and, in the case of the Bell distribution, the value of ${h^{-1}(\widehat\mu_n)}$ is obtained by using the function \texttt{uniroot} of the package Stats.
     
	\subsection{Empirical significance level}
\begin{figure}[H]
	\centering
	\includegraphics[width=0.9\textwidth]{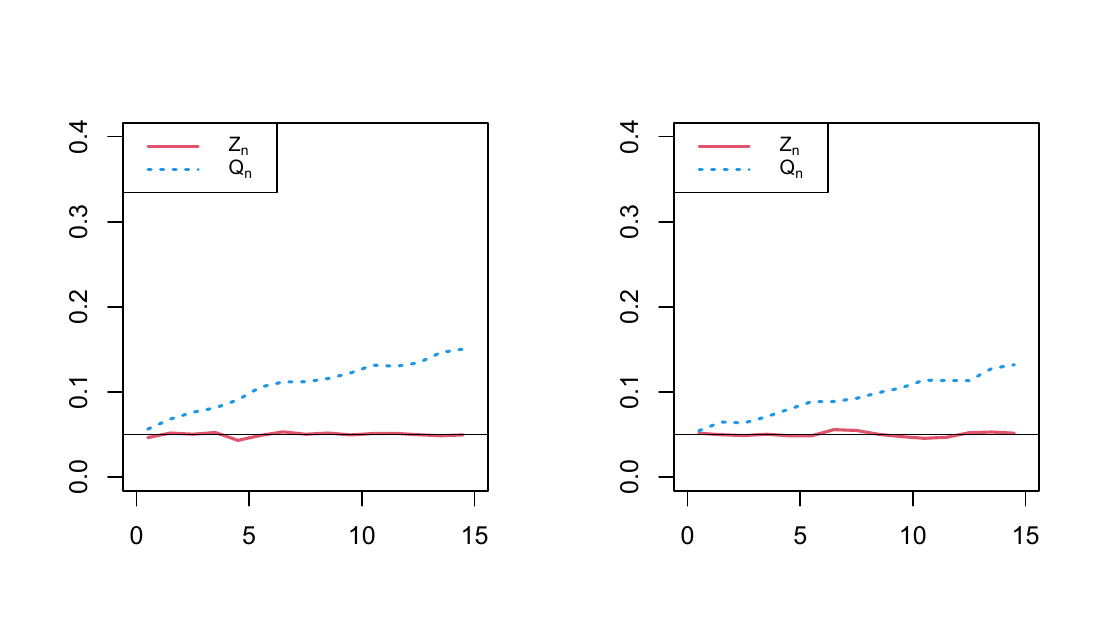}
	\vskip-1.4cm
	\caption{Empirical significance level under the shifted Borel distribution for $n=30$ (left panel) and $n=50$ (right panel). In the abscissa $\mu$ values.}
	\label{fig:figure1}
\end{figure}

\begin{figure}[H]
	\centering
	\includegraphics[width=0.9\textwidth]{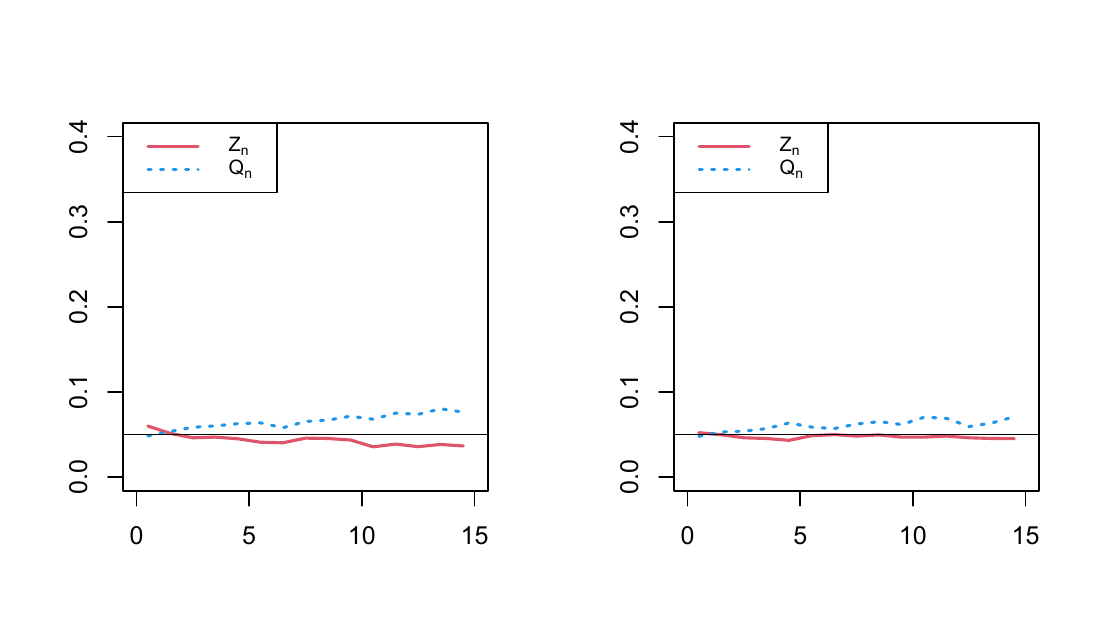}
	\vskip-1.4cm
	\caption{Empirical significance level under the Geometric distribution for $n=30$ (left panel) and $n=50$ (right panel). In the abscissa $\mu$ values.}
	\label{fig:figure2}
\end{figure}

\begin{figure}
	\centering
	\includegraphics[width=0.9\textwidth]{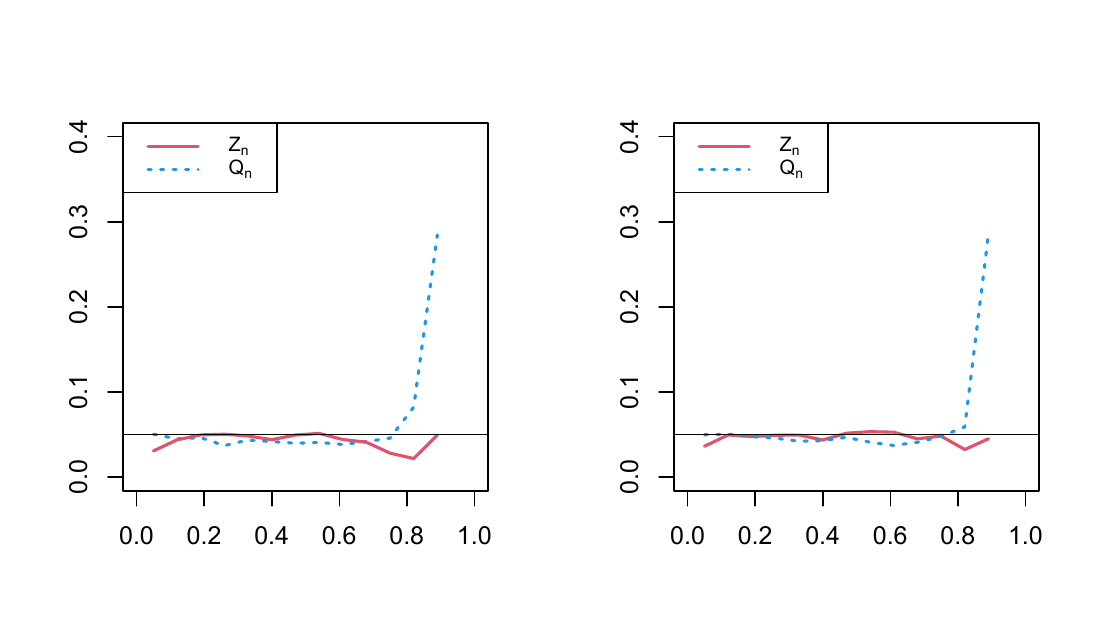}
		\vskip-1.4cm
	\caption{Empirical significance level under the Bell distribution for $n=30$ (left panel) and $n=50$ (right panel). In the abscissa values of $P(X=0)$.}
	\label{fig:figure3}
\end{figure}

	 First of all, we focus on empirically evaluating the actual significance level of the test. To this purpose, fixed the nominal level $\alpha=0.05$, $5000$ samples of size $ n=30,50$ are independently generated from the shifted Borel, Geometric and Bell distributions and the empirical significance level is computed as the proportion of rejections of the null hypothesis both for $Z_n$ and the chi-squared statistic $Q_n$. In particular, Figure~\ref{fig:figure1} and Figure~\ref{fig:figure2} display the empirical significance level for the shifted Borel and the Geometric distribution, respectively, for values of $\mu$ varying from $0.5$ to $15$ by $1$. Figure~\ref{fig:figure3} show the empirical significance level for the Bell distribution for probability of zero varying from $0.05$ to $0.89$ by $0.07$.

As to the shifted Borel (Figure \ref{fig:figure1}), the proposed test shows an empirical significance level almost equal to the nominal one for any $\mu$ already for the smaller sample size. On the other hand, the chi-squared test does not have a satisfactory behaviour even for the larger sample size, especially for large values of $\mu$, probably owing to the slow rate of convergence of $Q_n$. Considering the Geometric distribution (Figure \ref{fig:figure2}), for $n=30$, $Z_n$ shows conservativeness for larger $\mu$ values while the empirical level of $Q_n$ is greater than the nominal one and increases as $\mu$ increases. On the other hand, for $n=50$ the empirical level reached by $Z_n$ is almost indistinguishable from the nominal one and, even if the performance of $Q_n$ greatly improves, it is not completely satisfactory especially for the larger values of $\mu$. It is worth noting that under the Geometric distribution $Z_n$ has a very simple expression leading to a straightforward implementation of the test.  
	
    From Figure~\ref{fig:figure3}, it is at once apparent that, when considering the Bell distribution, the performances of both tests are quite satisfactory when $P(X=0)$ ranges between 0.1 and 0.8. Moreover, for $P(X=0)$ near to 0.9, while the empirical level of $Z_n$ remains rather close to the nominal one, especially for $n=50$, the empirical level of $Q_n$ dramatically increases. Moreover, the empirical level of both tests is very far from the nominal one for values of the probability close to 0 and close to 1.

	\subsection{Empirical power}
	The empirical power of the test based on $Z_n$ is investigated and compared to that of the chi-squared test considering, under the null hypothesis, the same distributions already adopted for assessing the empirical significance level. The power behaviour is assessed against some common alternative distributions with various parameter values and against contiguous alternatives, as introduced in Section \ref{sec:contig}. In order to highlight the behaviour of the proposed test in the most general context, we considered 
	alternative distributions for which consistency holds, but not necessarily belonging to the class 
	ensuring constant sign of the derivative of the ratio of the p.g.f.s.
	Alternative distributions include overdispersed and underdispersed, mixtures and zero-inflated distributions, together with distributions having mean close to variance. In particular, as well as in \cite{gurtler2000recent} and \cite{meintanis2008class}, we consider Poisson distribution denoted by $\mathcal{P}(\lambda)$, 
	Mixture of two Poisson  by $\mathcal{MP}(\lambda_{1},\lambda_{2})$ with mixture weight $0.5$, Binomial by $\mathcal{B}(k, p)$, Negative Binomial by $\mathcal{NB}(k, p)$, Generalized Hermite by $\mathcal{GH}(a, b, k)$, Discrete Uniform in $\{0,1,\dots,\nu\}$ by $\mathcal{DU}(\nu)$,
	Logarithmic Series by $\mathcal{LS}(\theta)$, Generalized Poisson by $\mathcal{GP}(\lambda_1,\lambda_2)$, Zero-inflated Binomial by $\mathcal{ZB}(k, p_{1}, p_{2})$, Zero-inflated Negative Binomial  by $\mathcal{ZNB}(k, p_{1}, p_{2})$, Zero-inflated Poisson by $\mathcal{ZP}(\lambda_{1},\lambda_{2})$,
	where various parameters values are considered (see Table~\ref{tab: tabella1} and ~\ref{tab: tabella2}).
	From each distribution, $5000$ samples of size $n=30, 50$ are independently generated and on each sample the tests based on $Z_n$ and $Q_n$ are performed at the nominal significance level $\alpha=0.05$. The empirical power of each test is computed as the proportion of rejections of the null hypothesis. Table~\ref{tab: tabella1} reports, for each model specified under null hypothesis and for each alternative distribution, the empirical powers achieved by the test based on $Z_n$ and $Q_n$, respectively, for $n=30$, while Table~\ref{tab: tabella2} for $n=50$.
	
	Not surprisingly, none of the two tests shows better performance with all models and all alternative distributions. Indeed, the power crucially depends both on the model specified under the null and alternative hypothesis and on the set of parameter values for alternatives in the same class. In particular, when under the null hypothesis the shifted Borel distribution is considered, the performance of both tests is rather satisfactory even with $n=30$, but the proposed test seems to be superior. If the Geometric distribution is specified under $H_0$, the power of both tests generally decreases and heavily deteriorates for zero-inflated alternative distributions, with the only remarkable exception for the chi-squared test when the zero-inflated binomial distribution is considered. Finally, a further decrease in the power of both tests occurs for the Bell distribution, with $Z_n$ reaching unbiasedness against the generalized Hermite alternatives only for $n=50$, even thought in this case both tests show a very poor performance.

\begin{figure}
	\centering
	\includegraphics[width=\textwidth]{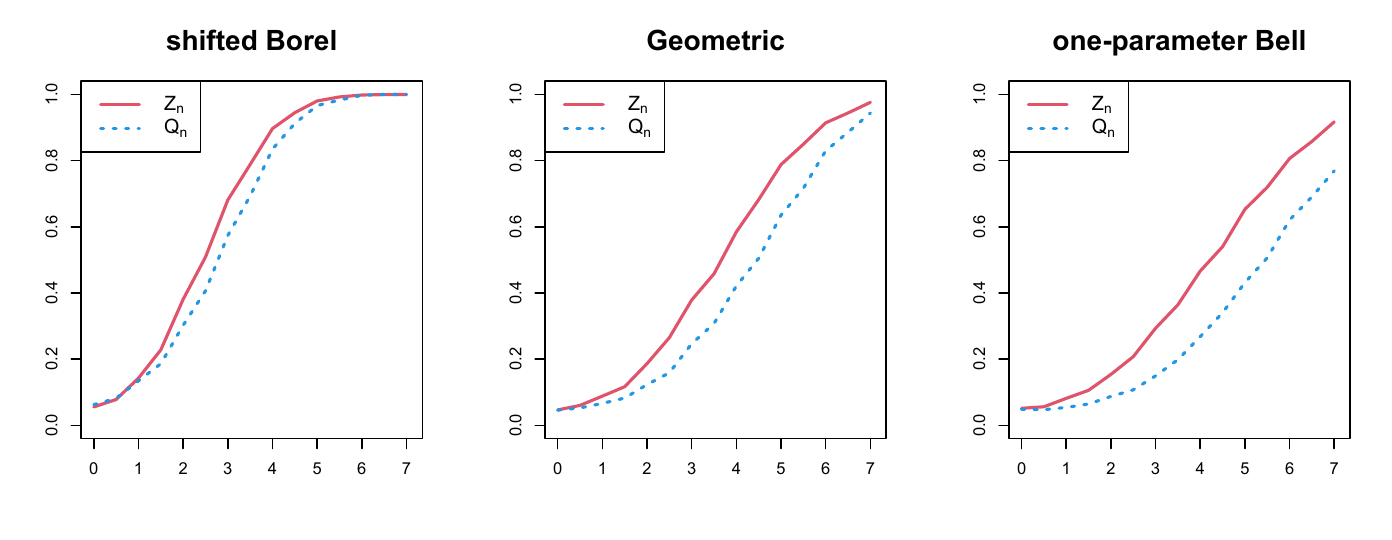}
	\vskip-0.85cm
	\caption{Empirical power under contiguous distributions for $n=50$. In the abscissa $\lambda$ values.}
	\label{fig:figure4}
\end{figure}

\begin{figure}
	\centering
	\includegraphics[width=\textwidth]{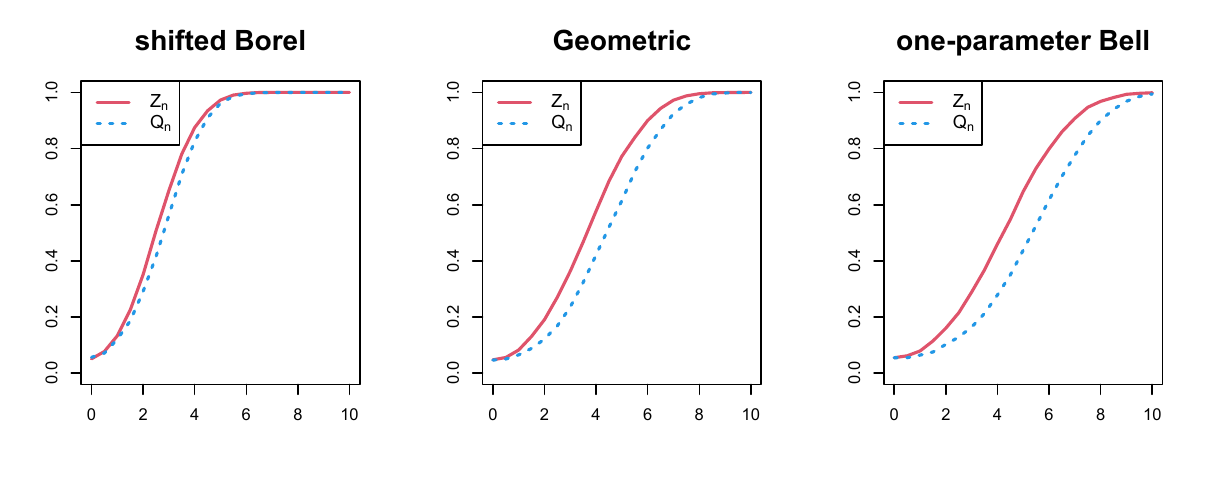}
	\vskip-1cm
	\caption{Empirical power under contiguous distributions for $n=100$. In the abscissa $\lambda$ values.}
	\label{fig:figure5}
\end{figure}

	As to the contiguous alternatives, for each distribution specified under the null hypothesis, shrinking mixtures are obtained according to \eqref{Xprime}. In particular, the component $X_l$ is a shifted Borel r.v., a Geometric r.v., a Bell r.v., all having $\mu=1$, respectively, while $Y_l\sim\mathcal{B}(4, 0.25)$ and $\lambda$ varies from $0$ to $\sqrt{n}$ by $0.5$. Figure \ref{fig:figure4} and Figure \ref{fig:figure5} show that the empirical power is rather satisfactory for both tests already for $n=50$ with a remarkable increase for $n=100$, with the best performance achieved with shifted Borel component. However, the empirical power of $Z_n$ is higher for any value of $\lambda$, any component and both sample sizes.

\begin{figure}
	\centering
	\includegraphics[width=0.9\textwidth]{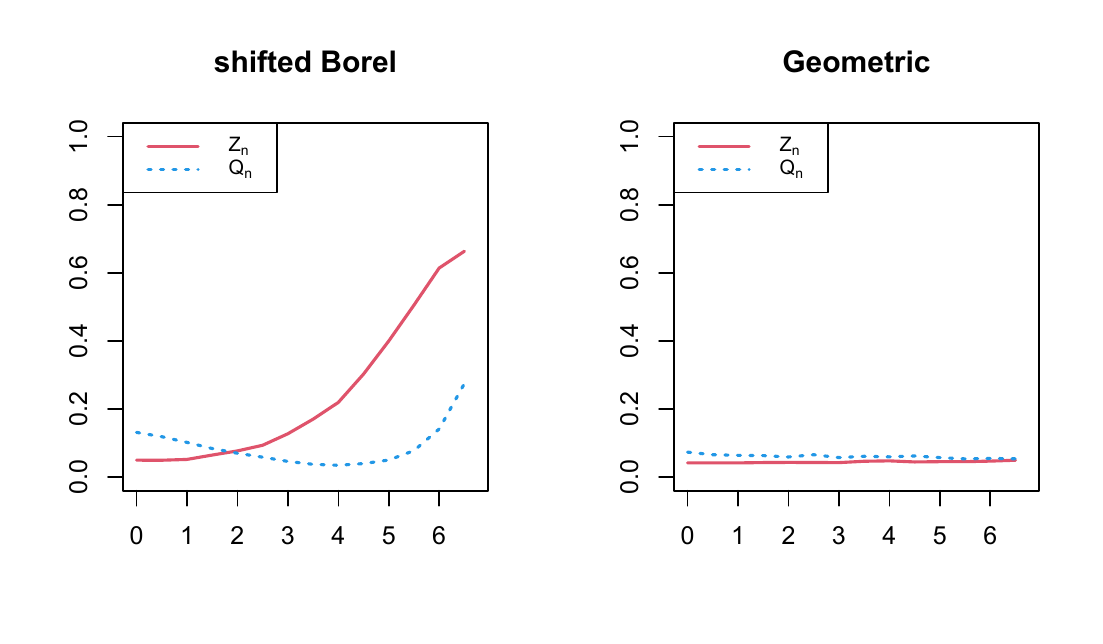}
	\vskip-1.4cm
	\caption{Empirical power under $\alpha$-fraction distributions for $n=50$. In the abscissa $\lambda$ values.}
	\label{fig:figure6}
\end{figure}

\begin{figure}
	\centering
	\includegraphics[width=0.9\textwidth]{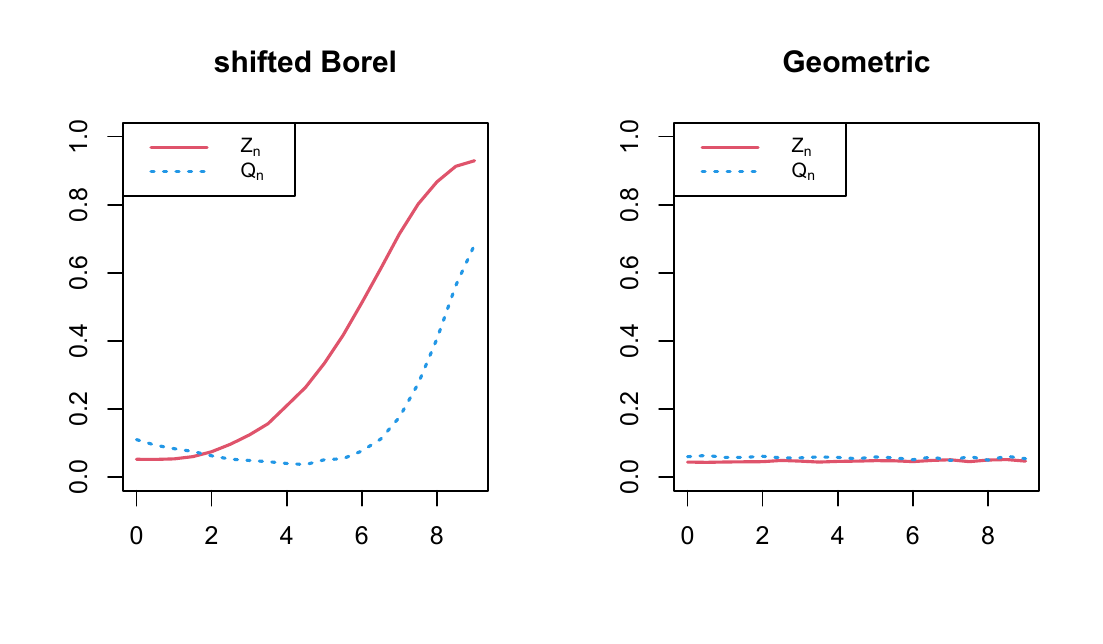}
		\vskip-1.4cm
	\caption{Empirical power under $\alpha$-fraction distributions for $n=100$. In the abscissa $\lambda$ values.}
	\label{fig:figure7}
\end{figure}

	Finally, Figure~\ref{fig:figure6} and Figure~\ref{fig:figure7} show the empirical power under alternatives obtained by means of the binomial thinning for $n=50$ and $n=100$. Since the thinning operator preserves the law in most of the cases, we report the Borel case in which the law is not preserved and the Geometric case in which the law is well-known to be preserved. In both cases for $n=50$, $\mu$ is fixed at $15$ and $\lambda$ varies from $0$ to $6.5$ by $0.5$ while, for $n=100,$ $\mu$ is fixed at $15$ and $\lambda$ varies from $0$ to $9$ by $0.5$. Coherently with the theoretical results, in the Geometric case the empirical power of both tests remains constantly close to the nominal level for both sample sizes. On the other hand, in the Borel case the empirical power of both tests start to increase when $\lambda$ increases but $Z_n$ performs better than $Q_n$ and shows a satisfactory behaviour for $n=100$.

Recently, \cite{jimenez2021test} proposed a test statistic $T_n$ for the Geometric distribution based on linear regression on order statistics which presents a competitive behaviour with respect to already existing tests. The bootstrap estimators of the null bootstrap distribution of $T_n$ are based on the probability density function of the normal law with zero mean and variance $\beta$. Indeed, their test statistic $T_n$ depends on the values of the parameter $\beta$. The authors recommend to take $\beta=1$ or $\beta=1.5$. Table~\ref{tab: tabella3} contains the empirical powers of $Z_n$ under the same alternative distributions considered in that paper, together with those reported by the authors for $T_n$ when $\beta=1$. Among the alternatives distributions, the Neyman type A distribution is denoted by $\mathcal{NA}(\lambda_1,\lambda_2)$ and the Discrete Weibull by $\mathcal{DW}(q,b)$. Notwithstanding $Z_n$ is compared to $T_n$, which is specifically tailored for the Geometric distribution, the empirical powers are rather similar for most of the alternatives. Exceptions are the Binomial distribution or the Zero-inflated Poisson distribution, where the performance of both tests more heavily depends on the parameters values.

\section{Discussion}
\label{sec:disc}
A huge literature deals with testing continuous distributions, while a few proposals have been developed for testing the fit of discrete distributions. Among these, many are tailored to deal with particular distributions and so they are of limited applicability, while there is an emerging need of tests allowing to specify an extremely broad class of distributions under the null hypothesis. Undoubtedly, the chi-squared test is the most widely adopted, notwithstanding the arbitrariness of its implementation due to the requirement on the frequency minimum values. Similarly to the chi-squared test, the proposed test shows considerable flexibility as it allows to specify any count distribution with finite second moment under the null hypothesis and can be considered deriving from a $L^1$ distance for the class of alternative distributions for which $D_\mu(s)$ has constant sign. The resulting test statistic has a simple expression, depending on the empirical p.g.f. only through the probability of zero occurrence, and an asymptotic normal distribution. Moreover, it must be pointed out that the test can be adopted for a much broader class of alternative distributions. Indeed, the test is consistent for the alternatives in the class but also for all the alternatives for which the probability of zero is different from that under the null hypothesis.   
The test performance are rather satisfactory even for moderate sample sizes, also compared to that of the test by \cite{jimenez2021test}, specifically tailored for the Geometric distribution. 
The test shows some criticalities when an inflation of zeroes occurs, which could be overcome introducing a family of test statistics, indexed by a parameter $h$, depending on the cumulative distribution function at $h$ instead of on the probability in zero. Further research will be devoted to this last issue and to the use of alternative estimators involved in the implementation of the statistic, in order to improve the performance of the test. Finally, the generalization of the test to families of count distributions indexed by a $k$-variate parameter will be investigated.

\bibliography{Riferimentibib_new}
	
\clearpage

\appendix
\section*{Appendix}
\section{Tables}

\begin{table}[h]
	\centering
	\caption{Empirical power (percent) with $5\%$ nominal significance level for $n=30$.}
	\medskip
	\renewcommand\arraystretch{1.2}
		\begin{tabular}{lcccccc} 
		\toprule
			& \multicolumn{6}{c}{Model under $H_0$} \\
		\cmidrule(lr){2-7}
			 & \multicolumn{2}{c}{shifted Borel}& \multicolumn{2}{c}{Geometric}&\multicolumn{2}{c}{Bell}\\	
		  	\cmidrule(lr){2-3} 		\cmidrule(lr){4-5} 		\cmidrule(lr){6-7} 						
			Alternative& $Z_{30}$ & $Q_{30}$& $Z_{30}$ &$Q_{30}$ & $Z_{30}$ & $Q_{30}$\\
			\midrule
			$\mathcal{P}(0.5)$           &   56.0&47.0&25.2&12.4&18.0&10.9          \\
			$\mathcal{P}(1)$             &   94.7&91.4&57.7&38.9&39.3&20.2        \\
			$\mathcal{P}(2)$             &   100.0&100.0&88.0&77.5&54.4&26.1           \\
			$\mathcal{MP}(1,2)$          &   99.0&97.8&62.4&43.4&31.7&14.3         \\
			$\mathcal{MP}(1,3)$          &   99.0&97.9&45.7&34.7&13.5&6.4           \\
			$\mathcal{MP}(1,4)$          &   98.3&96.8&26.6&26.6&5.3&6.7         \\
			$\mathcal{B}(4,0.25)$        &   99.4&98.9&85.9&72.4&73.6&50.4      \\
			$\mathcal{B}(30,0.1)$        &   100.0&100.0&97.8&97.2&60.2&34.4           \\
			$\mathcal{NB}(4,0.75)$       &   94.2&89.5&39.8&24.8&17.4&8.4            \\
			$\mathcal{NB}(10,0.9)$       &   94.7&90.7&51.2&32.1&30.2&13.5        \\
			$\mathcal{GH}(1,1.25,2)$     &   100.0&100.0&43.0&62.3&4.0&6.1             \\
			$\mathcal{GH}(1,1.5,2)$      &   100.0&100.0&44.7&69.8&3.3&7.2           \\
			$\mathcal{DU}(3)$            &   99.0&100.0&61.5&94.9&29.2&67.9          \\
			$\mathcal{LS}(0.6)$          &   100.0&100.0&100.0&100.0&100.0&100.0                \\
			$\mathcal{LS}(0.8)$          &   100.0&100.0&100.0&99.4&99.8&97.7         \\
			$\mathcal{GP}(1,0.1)$        &   90.4&84.4&38.7&23.5&21.5&9.9         \\
			$\mathcal{GP}(3,0.25)$       &   100.0&100.0&75.5&68.5&10.1&6.0              \\
			$\mathcal{ZB}(5,0.9,0.2)$    &   93.3&100.0&2.8&100.0&33.6&100.0          \\
			$\mathcal{ZNB}(5,0.9,0.1)$   &   36.7&30.4&11.7&6.6&7.9&5.5          \\
			$\mathcal{ZP}(1,0.2)$        &   65.2&60.3&18.3&12.7&9.8&5.4        \\
		\bottomrule
		\end{tabular}  

	\label{tab: tabella1}
\end{table}

\clearpage

\begin{table}
	\centering
	\caption{Empirical power (percent) with $5\%$ nominal significance level for $n=50$.}
	\medskip
	\renewcommand\arraystretch{1.2}
		\begin{tabular}{lcccccc} 
			\toprule
			& \multicolumn{6}{c}{Model under $H_0$} \\
			\cmidrule(lr){2-7} 
			 & \multicolumn{2}{c}{shifted Borel}& \multicolumn{2}{c}{Geometric}&\multicolumn{2}{c}{Bell}\\	
		 	\cmidrule(lr){2-3} 		\cmidrule(lr){4-5} 		\cmidrule(lr){6-7} 							
		Alternative	&  $Z_{50}$ & $Q_{50}$& $Z_{50}$ & $Q_{50}$ & $Z_{50}$ & $Q_{50}$\\
			\midrule
			$\mathcal{P}(0.5)$                  &    77.8&70.4&37.0&24.7&28.9&16.8             \\
			$\mathcal{P}(1)$                    &    99.6&99.0&79.7&64.2&60.4&34.9             \\
			$\mathcal{P}(2)$                    &    100.0&100.0&98.7&96.4&78.7&50.6             \\
			$\mathcal{MP}(1,2)$                 &    100.0&99.9&83.7&70.0&49.5&25.0                \\
			$\mathcal{MP}(1,3)$                 &    100.0&100.0&68.3&54.1&19.6&9.0              \\
			$\mathcal{MP}(1,4)$                 &    100.0&99.9&42.4&40.4&5.4&7.8             \\
			$\mathcal{B}(4,0.25)$               &    100.0&100.0&97.6&95.2&92.7&78.8             \\
			$\mathcal{B}(30,0.1)$               &    100.0&100.0&100.0&100.0&87.1&69.2                \\
			$\mathcal{NB}(4,0.75)$              &    99.7&98.4&58.9&41.2&26.5&12.7          \\
			$\mathcal{NB}(10,0.9)$              &    99.7&98.9&72.3&54.3&47.1&24.6          \\
			$\mathcal{GH}(1,1.25,2)$            &    100.0&100.0&67.9&85.6&5.4&8.7             \\
			$\mathcal{GH}(1,1.5,2)$             &    100.0&100.0&72.4&91.6&5.0&9.2            \\
			$\mathcal{DU}(3)$                   &    100.0&100.0&85.5&100.0&49.5&97.6             \\
			$\mathcal{LS}(0.6)$                 &    100.0&100.0&100.0&100.0&100.0&100.0                 \\
			$\mathcal{LS}(0.8)$                 &    100.0&100.0&100.0&100.0&100.0&100.0                 \\
			$\mathcal{GP}(1,0.1)$               &    98.7&96.7&58.7&39.2&31.7&14.7          \\
			$\mathcal{GP}(3,0.25)$              &    100.0&100.0&95.3&91.3&20.9&8.8                \\
			$\mathcal{ZB}(5,0.9,0.2)$           &    99.2&100.0&3.4&100.0&53.0&100.0              \\
			$\mathcal{ZNB}(5,0.9,0.1)$          &    53.3&48.3&15.7&10.1&10.9&6.7          \\
			$\mathcal{ZP}(1,0.2)$               &    86.0&81.2&26.7&20.3&14.5&7.8             \\
			\bottomrule
		\end{tabular}  
	\label{tab: tabella2}
\end{table}	

\clearpage

\begin{table}
	\centering
	\caption{Empirical power (percent) with $5\%$ nominal significance level of tests for the Geometric distribution for $n=20$ and $n=40$.}
	\medskip
	\renewcommand\arraystretch{1.2}
		\begin{tabular}{lcccc} 
			\toprule
			 & \multicolumn{4}{c}{Geometric}\\	
			\cmidrule(lr){2-5}
			
		Alternative	& $Z_{20}$ & $T_{20}$& $Z_{40}$ &$T_{40}$\\
		\midrule
			
			$\mathcal{P}(0.5)$           & 14  & 13 & 29  & 24   \\
			$\mathcal{P}(1)$             & 42  &  37&  70 &  68  \\
			$\mathcal{P}(1.5)$           &  60 & 63 & 88  &  91  \\
			$\mathcal{ZP}(2,0.8)$        & 44  & 21 & 83  & 44   \\
			$\mathcal{ZP}(1,0.1)$           &  25 & 23 & 44  & 43   \\
			$\mathcal{ZP}(2,0.2)$           & 14  & 26 & 24  &  48  \\
			$\mathcal{NA}(5,0.1)$           & 9  & 9 & 17  &  15  \\
			$\mathcal{NA}(5,0.2)$           & 25  & 22 & 44  & 42   \\
			$\mathcal{NA}(5,0.3)$           &  35 & 35 & 59  &  65  \\
			
			$\mathcal{DW}(1.4,0.4)$           & 16  & 15 & 32  & 28   \\
			$\mathcal{DW}(1.4,0.6)$           &  27 &23  & 48  & 46   \\
			$\mathcal{DW}(1.4,0.8)$           &  28 &  32&  51 & 60   \\
			
			$\mathcal{B}(10,0.05)$           & 17  & 16 & 36  & 29   \\
			$\mathcal{B}(10,0.1)$            & 14  & 46 &  33 & 79   \\
			$\mathcal{B}(10,0.15)$           &  74 & 75 & 96  & 98   \\
			
			$\mathcal{NB}(5,0.89)$         & 13  & 12 & 26  & 22   \\
			$\mathcal{NB}(5,0.83)$           &  26 & 24 &  46 & 44   \\
			$\mathcal{NB}(5,0.75)$           & 38  & 38 & 67  & 71   \\
			
		\bottomrule
		\end{tabular}  

	\label{tab: tabella3}
\end{table}

\clearpage

\section{Figures}

\end{document}